\documentclass{au}
\usepackage{amssymb,amsfonts,amscd}
\theoremstyle{plain}
\newtheorem{thm}{Theorem}[section]

\newtheorem{prop}[thm]{Proposition}

\theoremstyle{definition}
\newtheorem*{ex}{Example}

\newtheorem{defi}[thm]{Definition}

\newtheorem*{rem}{Remark}

\begin{document}

\author{E. Graczynska}
\address{Technical University of Opole, Institute of Mathematics
\newline
ul. Luboszycka 3, 45-036 Opole, POLAND}
\email{egracz@po.opole.pl \hspace{1cm}
http://po.opole.pl/$\sim$egracz/}
\author{D. Schweigert}
\address{Universit\"{a}t Kaiserslautern, Fachbereich Mathematik
\newline
Postfach 3049 \\ 67653 Kaiserslautern, Germany}
\email{schweige@mathematik.uni-kl.de}

\title{Hyperquasivarieties}

\begin{abstract}
We consider the notion of hyper-quasi-identities and hyperquasivarieties,
as a common generalization of the concept of quasi-identity and quasivariety invented by
A. I. Mal'cev, cf. \cite{AIM}, cf. \cite{VAG} and hypervariety invented by the
authors in \cite{4}.


\end{abstract}

\maketitle

\section{Notations}

An identity is a pair of terms where the variables are bound by
universal quantifiers. Let us take the following medial identity
as an example
\[
\forall u\forall x\forall y\forall w\,(u\cdot x)\cdot (y\cdot
w)=(u\cdot y)\cdot (x\cdot w).
\]
Let us look at the following hyperidentity
\[
\,\forall F\forall u\forall x\forall y\forall
w\,F(F(u,v),F(x,y))=F(F(u,x),F(v,y)).
\]
The hypervariable $F$ is considered in a very specific way.
Firstly every hypervariable is restricted to functions of a given
arity. Secondly $F$ is restricted to term functions of the given type.
Let us take the variety $A_{n,0}$ of abelian groups of finite exponent
$n.$ Every binary term $t\equiv t(x,y)$ can be presented by
$t(x,y)=ax+by$ with $a,b\in \mathbb N_{0}.$ If we substitute the
binary hypervariable $F$ in the above hyperidentity by $ax+by,$
leaving its variables unchanged, we get
\[
a(au+bv)+b(ax+by)=a(au+bx)+b(av+by).
\]
This identity holds for every term $t(x,y)=ax+by$ for the variety
$A_{n,0}$. Therefore we say that the hyperidentity holds for the
variety $A_{n,0}.$

\section{Hyper-quasi-identities}

In the sequel we use the definition of {\em hyperterm} from
\cite{4}. We accept the notation of \cite{PMC}, \cite{5},
\cite{VAG} and \cite{DHRK}.

We recall only our definitions of \cite{4} of the fact that a hyperidentity
is satisfied in an algebra of a given type and the notion of hypervariety:

\begin{defi}\label{D:2.1}
An algebra ${\bf A}$ satisfies a hyperidentity $h_{1}=h_{2}$ if
for every substitution of the hypervariables by terms (of the same
arity) of ${\bf A}$ leaving the variables unchanged, the
identities which arise hold in ${\bf A}$. In this case, we write
${\bf A} \models (h_{1}=h_{2})$. A variety $V$ satisfies a hyperidentity
$h_{1}=h_{2}$ if every algebra in the variety does;
in this case, we write $V\models (h_{1}=h_{2})$.
\end{defi}
\begin{defi}\label{D:2.2}
A class $V$ of a algebras of a given type is called a hypervariety if
and only if $V$ is defined by a set of hyperidentities.
\end{defi}
The following was proved in \cite{4}:

\begin{thm} \label{T:2.3}
A variety $V$ of type $\tau$ is defined by a set of hyperidentities
if and only if $V = HSPD(V)$, i.e. $V$ is a variety closed under derived
algebras of type $\tau$.
\end{thm}
We recall from \cite{AIM} and \cite{VAG}:
\begin{defi}\label{D:2.4}
A quasi-identity $e$ is an implication of the form:
\begin{center}
$(t_{0} = s_{0}) \wedge ... \wedge (t_{n-1} = s_{n-1}) \rightarrow (t_{n} = s_{n})$.
\end{center}
where $t_{i} = s_{i}$ are $k$-ary identities of a given type, for $i = 0,...,n$.

A quasi-identity above is {\em satisfied in an algebra} ${\bf A}$ of a given type
if and only if the following implication is satisfied in ${\bf A}$:  given a 
sequence $a_{1},...,a_{k}$ of elements of $A$. If this elements
satisfy the equations $t_{i}(a_{1},...,a_{k}) = s_{i}(a_{1},...,a_{n})$
in ${\bf A}$, for $i = 0, 1,..., n-1$, then the equality
$t_{n}(a_{1},...,a_{k}) = s_{n}(a_{1},...,a_{k})$ is satisfied in ${\bf A}$.
In that case we write:
\begin{center}
${\bf A} \models (t_{0} = s_{0}) \wedge ... \wedge (t_{n-1} = s_{n-1}) \rightarrow (t_{n} = s_{n})$.
\end{center}

A quasi-identity $e$ is {\em satisfied in a clas}s $V$ of algebras of a given type,
if and only if it is satisfied in all algebras ${\bf A}$ belonging to $V$.

\end{defi}

We modify the definition above in the following way:
\begin{defi}\label{D:2.5}
A hyper-quasi-identity $e^{*}$ is an implication of the form:
\begin{center}
$(T_{0} = S_{0}) \wedge ... \wedge (T_{n-1} = S_{n-1}) \rightarrow (T_{n} = S_{n})$.
\end{center}
where $T_{i} = S_{i}$ are hyperidentities of a given type, for $i = 0,...,n$.
                           
A hyper-quasi-identity $e^{*}$ is hyper-satisfied (holds) in an algebra ${\bf A}$ is
and only if the following implication is satisfied:\\
If $\sigma$ is a hypersubstitution of type $\tau$ and the elements
$a_{1},...,a_{n} \in A$ satisfy the equalities
$\sigma(T_{i})(a_{1},...,a_{k}) = \sigma(S_{i})(a_{1},...,a_{k})$ in ${\bf A}$,
for $n = 0,1,...,n-1$, then the equality
$\sigma(T_{n})(a_{1},...,a-{k}) = \sigma(S_{n})(a_{1},...,a_{k})$
holds in ${\bf A}$.
\end{defi}
By other words, hyper-quasi-identity is a universally closed  Horn
$\forall x \forall \sigma$-formulas, where x vary over all sequences of
individual variables (occuring in terms of the implication) and $\sigma$
vary over all hypersubstitutions of a given type.

\begin{rem}
All hyper-quasi-identities and hyperidentities are written without quantifiers
but they are considered as universally closed Horn $\forall$-formulas
(cf. \cite{AIM}). A syntactic side of the notions described here will be
considered in a forthcoming paper.
\end{rem}

\section{Examples of hyper-quasi-identities}
\subsection{Hyper-quasi-identities for abelian algebras}
\begin{defi}\label{D:3.1}
An algebra ${\bf A}$ is called abelian if for every $n > 1$ and every n-ary term
operation $f$ of ${\bf A}$ and for all $u,v,x_{1},...,x_{n-1},y_{1},...,y_{n-1}$
the following equivalence holds:
\begin{center}
$f(u,x_{1},...,x_{n}) = f(u,y_{1},...,y_{n}) \leftrightarrow f(v,x_{1},...,x_{n}) = f(v,x_{1},...,x_{n})$
\end{center}
A variety $V$ is called abelian, if each algebra of $V$ is abelian.
\end{defi}
It follows from  \cite[p. 40]{DHRK}, \cite[p. 290]{2}:

\begin{prop}\label{P:3.1}
An algebra ${\bf A}$ is abelian if and only if the following  hyper-quasi-identity
holds in ${\bf A}$:
\begin{center}
$F(u,x_{1},...,x_{n}) = F(u,y_{1},...,y_{n}) \leftrightarrow
F(v,x_{1},...,x_{n}) = F(v,y_{1},...,y_{n})$.
\end{center}
\end{prop}
\begin{ex}
The variety $RB$ of rectangular bands fulfills the hyperidentities of type
$(1,2,3,...,n,...)$:\\
$F(x,...,x) = x$, $F(F(x_{11},...,x_{1n}),x_{2},...,x_{n}) = F(x_{11},x_{2},...,x_{n})$, \\
$F(F(x_{1},...,x_{n-1}),F(x_{n1},...,x_{nn})) = F(x_{1},...,x_{n-1},x_{nn})$.\\
We can derive by the associative hyperidentity:
$F(x,F(y,z)) = (F(x,y),z)$, that the following hyper-quasi-identity holds in $RB$:\\
$F(u,x_{1},...,x_{n}) = F(u,y_{1},...,y_{n}) \rightarrow
F(v,u,x_{1},...,x_{n-1}) = F(v,u,y_{1},...,y_{n-1}) \rightarrow
F(v,x_{1},...,x_{n}) = F(v,y_{1},...,y_{n})$, i.e. the variety  $RB$ is abelian.
\begin{rem}
For further examples of abelian varieties consult \cite{DHRK} and \cite{11}.
\end{rem}

\end{ex}
\subsection{Semidistributive lattices}
A lattice is {\em joinsemidistributive} if it satisfies the following
condition, cf. \cite[p. 82]{DHRK},  \cite[p. 141]{VAG}:
\begin{center}
($SD_{\vee}$) $x \vee y = x \vee z$ implies $x \vee y = x \vee (y \wedge z)$
\end{center}
The {\em meetsemidistributivity} is defined by duality:
\begin{center}
($SD_{\wedge}$) $x \wedge y = x \wedge z$ implies $x \wedge y = x \wedge (y \vee z)$
\end{center}

A lattice is {\em semidistributive} if it is simultaneously join and meet
semidistributive.

\begin{prop}\label{P:3.2}
Let ${\bf L} = (L, \wedge, \vee)$ be a lattice which is semidistributive.
Then the following hyper quasi-identity is hypersatisfied in ${\bf L}$:
\begin{center}
$(F(x,y) = F(x,z)) \rightarrow (F(x,y) = F(x,G(y,z)))$
\end{center}
\end{prop}
\begin{proof}
We consider all cases on hypervariables of a semidistributive lattice.

Case 1. $F(x,y) := x$.\\
Obviously $(x=x) \rightarrow (x =x)$.

Case 2. $F(x,y) := y$. Consider $G(y,z) := y, z, y \wedge z, y \vee z$.\\
Then the following quasi-identities are satisfied in ${\bf L}$: \\
$(y = z) \rightarrow (y = y)$, $(y =z) \rightarrow (y = z)$,
$(y = z) \rightarrow (y = y \wedge z)$,\\
$(y = z) \rightarrow (y = y \vee z)$.

Case 3. $F(x,y) := x \wedge y$. Consider $G(y,z) := y,z, y\wedge z, y \vee z$. \\
Then the following quasi-identities hold in ${\bf L}$: \\
$(x \wedge y = x \wedge z) \rightarrow  (x \wedge y = x \wedge y)$,
$(x \wedge y = x \wedge z) \rightarrow (x \wedge y = x \wedge z)$, \\
$(x \wedge y = x \wedge z) \rightarrow (x \wedge y = x \wedge (y\wedge z))$,
$(x \wedge y = x \wedge z) \rightarrow (x \wedge y = x \wedge (y \vee z))$.
The last one follows from the  meet semidistributivity of ${\bf L}$.

Case 4.  $F(x,y) := x \vee y$. Consider $G(y,z) := y,z, y\wedge z, y \vee z$.\\
Then the following quasi-identities hold in ${\bf L}$:\\
$(x \vee y = x \vee z) \rightarrow  (x \vee y = x \vee y)$,
$(x \vee y = x \vee z) \rightarrow (x \vee y = x \vee z)$, \\
$(x \vee y = x \vee z) \rightarrow (x \vee y = x \vee (y\wedge z))$,
$(x \vee y = x \vee z) \rightarrow (x \vee y = x \vee (y \vee z))$.
The prelast one follows from the  join semidistributivity of ${\bf L}$.
\end{proof}

\section{Derived algebras}
We recall from \cite{4}, cf. \cite[p. 145]{PMC} the notion of a
{\em derived algebra} and the {\em derived class} of algebras.
Given a type $\tau = (n_{0},n_{1},...,n_{\gamma},...)$.
An algebra ${\bf B}$ of type $\tau$ is called a derived algebra of
${\bf A} = (A,f_{0},f_{1},...,f_{\gamma},...)$ if there exist term
operations $t_{0},t_{1},...,t_{\gamma},...$ of type $\tau$ such that
${\bf B} = (A,t_{0},t_{1},...,t_{\gamma},...)$. For a class $K$ of algebras
of type $\tau$ we denote by ${\bf D}(K)$ the class of all derived algebras
(of type $\tau$) of $K$. \\
A class $K$ is called {\em deriverably closed} if and only if ${\bf D}(K) \subseteq K$.

For a given algebra ${\bf A}$ we denote by $QId({\bf A})$ and $HQId({\bf A})$
the set of all quasi-identities and hyper-quasi-identities satisfied
(hypersatisfied) in ${\bf A}$, respectively.  Similarly for a class $K$, $QId(K)$ and $HQId(K)$
denote the set of all quasi-identities and hyper-quasi-identities satisfied
(hypersatisfied) in $K$, respectively. Following \cite{VAG} by ${\bf Q}(K)$
and $H{\bf Q}(K)$ we denote the class of all algebras of a given type satisfying (hypersatisfying) all the
quasi-identities and hyper-quasi-identities of $K$, respectively. The transfomation
$T$ of an quasi-identity $e$ into hyper-quasi-identity $e^{*}$ is defined in a natural way.
Similarly $T^{-1}$ transform every hyper-quasi-identity $e^{*}$ into the qasiidentity $e$.
(cf. \cite{4}).

\begin{prop}\label{P:4.1}
Given a class $K$ of algebras of type $\tau$. Then the following equality
holds:
\begin{center}
$HQId(K) = T(QId({\bf D}(K)))$.
\end{center}
\end{prop}
\begin{proof}
To prove the inclusion $\subseteq$, given a hyper-quasi-identity $e^{*}$
of $K$ and an algebra a derived algebra ${\bf B} = {\bf A}^{\sigma}$ for
${\bf A} \in K$. Then by definition 2.5 $T^{-1}(e)$ is satisfied in ${\bf B}$
as a quasi-identity, i.e. $e^{*} \in T(QId({\bf D}(K)))$.

For a proof of the converse inclusion, given a quasi-identity $e$
satisfied in ${\bf D}(K)$ and an algebra ${\bf A}$ of $K$. Let $a_{1},...,a_{k} \in A$.
Given an hypersubstitution $\sigma$ of type $\tau$ and consider $\sigma(e)$
in ${\bf A}$. As ${\bf A} \in {\bf D}(K)$, then $\sigma(e)$ may be considered as
a quasi-identity of the derived algebra ${\bf A}^{\sigma}$ (cf. \cite{10}).
Assume that $p_{i}^{{\bf A}^{\sigma}}(a_{1},,,a_{k}) = q_{i}^{{\bf A}^{\sigma}}(a_{1},...,a_{k})$, therefore
as ${\bf A} \in {\bf D}(K)$ we obtain that $p_{n}^{{\bf A}^{\sigma}}(a_{1},...,a_{k}) = q_{n}^{{\bf A}^{\sigma}}(a_{1},...,a_{n})$,
i.e. that the $T(e)$ holds in $K$ as a hyper-quasi-identity, i.e. $T(e) \in HQId(K)$.
\end{proof}
The role of derived algebras in solvability questions may be visualized by the
following:
\begin{thm}\label{T:4.2}
Given a locally finite variety $V$. The class $K$ of all locally solvable
algebras in $V$ is a hypervariety.
\end{thm}
\begin{proof}
By corollary 7.6, of \cite{DHRK} the class of all locally solvable algebras in a locally
finite variety $V$ is a variety. Due to our theorem 2.3 we need to show that
$K$ is deriverably closed. Given a derived algebra ${\bf B}$ of an algebra
${\bf A} \in K$ and its finite subalgebra ${\bf F}^{*}$. Then ${\bf F}^{*}$
is a derived subalgebra of a finite subalgebra ${\bf F}$ of ${\bf A}$, i.e.
${\bf F}^{*} = {\bf F}^{\sigma}$, for a hypersubstitution $\sigma$.
From theorem 5.7 of \cite{10} we conclude, that as ${\bf F}$ is solvable,
therefore ${\bf F}^{\sigma}$ is solvable.
\end{proof}
\begin{rem}
The assumption that $V$ is locally finite is essential, as for the variety of
abelian groups $AG$, which are locally solvable, it is easy to notice that
$AG$ does not constitute a hypervariety, as the abelian law: $xy = yx$ is
not satisfied in $AG$ as a hyperidentity.
\end{rem}
\subsection{Quasicompact classes}


We accept the definition of a {\em quasicompact classes} invented by
V. A. Gorbunov  \cite[p. 77]{VAG} and his theorem 2.3.1, which states:
\begin{thm}\label{T:4.3}
For a class $K$,  the prevariety ${\bf S}{\bf P}(K)$ is a quasivariety if and only if
$K$ is quasicompact.
\end{thm}
We reformulate the definition of \emph{ prevariety} to a \emph{hyperprevariety}
and  theorem above for the case of hyper-quasi-identities:
\begin{defi}\label{D:4.4}
For a a class $K$ of algebras, the class ${\bf S}{\bf P}({\bf D}(K))$
(for short ${\bf S}{\bf P}{\bf D}(K)$)
will be called hyperprevariety generated by $K$.
\end{defi}
\begin{def}\label{D:4.5}
A class $K$ of algebras of a given type is hyper-quasi-compact provided
that from the infinite implication, indexed by $I$ and hyper-satisfied in $K$
it follows that for some finite subset $F \subseteq I$ the finite implication
(restricted to $F$) is hypersatisfied in $K$.
\end{def}
\begin{thm}\label{T:4.6}
For a class $K$, the prevariety ${\bf S}{\bf P}({\bf D}(K))$ is a quasivariety if and only
if ${\bf D}(K)$ is quasicompact if and only if $K$ is hyper-quasi-compact.
\end{thm}
\begin{proof}
The first equivalence is the V. A. Gorbunov theorem of \cite[p. 77]{VAG}. \\
The second equivalence follows from the proposition 4.1.
\end{proof}
\section{Hyper-quasi-varieties}
We reformulate the notion of {\em quasivariety} invented by A. I. Mal'cev
in \cite[p. 210]{AIM} for the case of {\em hyper-quasi-identities} of
a given type in a natural way:
\begin{defi}\label{D:5.1}
A class $K$ of algebras of type $\tau$ is called  a hyperquasivariety if there
is a set  $\Sigma$ of hyper-quasi-identities of type $\tau$ such that
$K$ conisits exactly of those algebras of type $\tau$ that hypersatisfy all
the hyper-quasi-identities of $\Sigma$.
\end{defi}
From proposition 4.1 we obtain immediately:
\begin{thm}\label{T:5.2}
A quasivariety $K$ of algebras given type is a hyperquasivariety if and only if
it is deriverably closed.
\end{thm}

In the sequel we use the standard notation: ${\bf S}$ for the operator of
creating \emph{subalgebras}, ${\bf P}$, ${\bf P_{s}}$, ${\bf P_{r}}$, ${\bf P_{u}}$
and ${\bf P{\omega}}$ for \emph{products}, \emph{subdirect products},
\emph{reduced products}, \emph{ultraproducts} and \emph{direct products of finite families of
structures}, respectively.  ${\bf L}$ and ${\bf L_{s}}$ will denote the operators of
\emph{direct limits} and \emph{superdirect limits}.

Following  A. I. Mal'cev \cite[p. 153, 215]{AIM} for a given class $K$
of algebras of a given type, by ${\bf S}(K)$, ${\bf P}(K)$, ${\bf P_{s}}(K)$, ${\bf P_{r}}(K)$,
${\bf P_{u}}(K)$ we denote the class of algebras isomorphic to all possible subalgebras, direct products,
subdirect products, reduced (filtered) products or ultraproducts of algebras
of $K$, respectively. ${\bf P_{\omega}}$ is the class of algebras
isomorphic to the direct products of finite families of structures of $K$.
Similarly, ${\bf L}(K)$ and ${\bf L_{s}}(K)$ will be the class of algebras isomorphic to
direct and superdirectlimits of algebras of $K$, respectively (cf. \cite[p.21]{VAG}.

Adding a trivial system to $K$ we obtain the class $K_{0}$.

We reformulate the resut of \cite{4}:
\begin{prop}\label{P:5.3}
Given a class $K$ of algebras, then the following inclusions holds:\\
1) ${\bf D}{\bf S}(K) \subseteq {\bf S}{\bf D}(K)$; \\
2) ${\bf D}{\bf P}(K) \subseteq {\bf P}{\bf D}(K)$;\\
3) ${\bf D}{\bf P_{\omega}}(K) \subseteq {\bf P_{\omega}}{\bf D}(K)$;\\
4) ${\bf D}{\bf P_{s}}(K) \subseteq {\bf P_{s}}{\bf D}(K)$;\\
5) ${\bf D}{\bf P_{r}}(K) \subseteq {\bf P_{r}}{\bf D}(K)$;\\
6) ${\bf D}{\bf P_{u}}(K) \subseteq {\bf P_{u}}{\bf D}(K)$;\\
7) ${\bf D}{\bf L}(K) \subseteq {\bf L}{\bf D}(K)$;\\
8) ${\bf D}{\bf L_{s}}(K) \subseteq {\bf L_{s}}{\bf D}(K)$.
\end{prop}
\begin{proof}
Obviously an isomorphism respects the inclusions above.  \\
The inclusions 1) and 2) were proved in \cite{4}.\\ 
3) and 4) immediately follows from 2).

To show 5) assume that an algebra ${\bf B} \in {\bf D}{\bf P_{r}}(K)$.
Let ${\bf B} = (A, t^{\bf A}_{0},t^{\bf A}_{1},...,t^{\bf A}_{\gamma},...)$
is a derived algebra of a reduced product ${\bf A} = (A, f^{\bf A}_{0},f^{\bf A}_{1},...,f^{\bf A}_{\gamma},...)$,
where $A = (\Pi_{i \in I}A_{i})/\sim_{F}$, for a set $I$, a family
$({\bf A}_{i})_{i \in I}$ of algebras of type $\tau$ from $K$ and a filter $F$ over $I$,
i.e. ${\bf A} = (\Pi_{i \in I}{\bf A}_{i})/F$, where for any $n$-ary functional
symbol $f$ the following holds  (cf. \cite[p. 13]{VAG}):
\begin{center}
$f^{\bf A}(a_{0}/F,...,a_{n-1}/F) = a_{n}/F$ if and only if
$\{i \in I: f^{{\bf A}_{i}}(a_{0},...,a_{n-1}) = a_{n} \} \in F$.
\end{center}
Note that by the induction on the complexity of a term $t$ of type $\tau$
one may show, that the following holds for any $n$-ary polynomial $t$ of type
$\tau$:
\begin{center}
$t^{\bf A}(a_{0}/F,...,a_{n-1}/F) = a_{n}/F$ if and only if
$\{i \in I: t^{{\bf A}_{i}}(a_{0},...,a_{n-1}) = a_{n} \} \in F$.
\end{center}
From the above equality it follows that the algebra
${\bf B} = \Pi_{i \in I}{\bf A}_{i}^{\sigma}/F$, where
${\bf A}_{i}^{\sigma} = (A_{i}, t_{0}^{{\bf A}_{i}},t_{1}^{{\bf A}_{i}},...,t_{\gamma}^{{\bf A}_{i}},...)$,
for $i \in I$ is a derived algebra of ${\bf A}_{i}$, i.e. ${\bf B} \in {\bf P_{r}}{\bf D}(K)$.

The inclusion 6) is an immediate consequence of 5).

To prove 7), let ${\bf B} \in {\bf D}{\bf L}(K)$, i.e. given a derived algebra
${\bf B}$ of an algebra ${\bf A}$, where ${\bf A} = (A,f_{0}^{\bf A},f_{1}^{\bf A},...,f_{\gamma}^{\bf A},...)$ is a direct limit of
a triple $\Lambda = (I, {\bf A}_{i},g_{ij})$, where $I = (I, \leq)$ is an
up-directed set, $({\bf A}_{i})_{i \in I}$ is a family of algebras of a
given type $\tau$ and $\{g_{ij} : i,j \in I , i \leq j\}$ is a family of
homomorphisms of ${\bf A}_{i}$ into ${\bf A}_{j}$ called a \emph{direct
spectrum} over $(I, \leq)$, cf. \cite[p. 17]{VAG}. For a given direct spectrum
$\Lambda = ({\bf A}_{i},g_{ij})$, we consider the quotient set
$A = \bigcup _{i \in I} A_{i} \times \{i\}/\equiv$, where
\begin{center}
$(a,i) \equiv (b,j)$ if and only if $(\exists k \in I)(i,j \leq k, g_{ik}(a) = g_{jk}(b))$.
\end{center}
Let $<a,i>$ denotes the equivalence class by $\equiv$ containing $(a,i)$.
The operations on $A$ are defined by setting for any operation symbol of type $\tau$:
\begin{center}
$f^{\bf A}(<a_{0},i_{0}>,...,<a_{n-1},i_{n-1}>) =
<f^{{\bf A}_{j}}(g_{i_{0}j}(a_{0}),...,g_{i_{n-1}j}(a_{n-1})),j>$.
\end{center}
Note, that as $g_{ij}$ are homomorphisms for all $i \leq j$, $i,j \in I$,
therefore for any polynomial symbol $p$ of type $\tau$, the following holds:
\begin{center}
$p^{\bf A}(<a_{0},i_{0}>,...,<a_{n-1},i_{n-1}>) =
<p^{{\bf A}_{j}}(g_{i_{0}j}(a_{0}),...,g_{i_{n-1}j}(a_{n-1})),j>$.
\end{center}
As ${\bf B} = {\bf A}^{\sigma}$ for a hypersubstitution $\sigma$ of type $\tau$,
therefore ${\bf B} \in {\bf L}{\bf D}(K)$, namely ${\bf B}$ is a direct limit of the
triple $\Lambda^{\sigma} = (I,{\bf A}_{i}^{\sigma},g_{ij})$ of derived algebras of
$K$, i.e. of a derived spectrum of $\Lambda$.

To prove 8), recall \cite[17]{VAG}, that a direct spectrum
$\Lambda = (I,{\bf A}_{i},gij)$ is called \emph{superdirect} if the mappings
$g_{ij}: A_{i} \rightarrow A_{j}$ are surjective. The family
$({\bf A}_{i})_{i \in I}$ is referred to as \emph{superdirect family}.
The direct  limit of a superdirect spectrum is called the \emph{superdirect limit}.
Note, that the derived spectrum $\Lambda^{\sigma}$ of  a superdirect spectrum
$\Lambda$ is superdirect. This together with 7) proves 8).
\end{proof}
Via corollary 2.3.4 \cite[p. 79]{VAG} and proposition 5.1 we obtain another
characterization of hyperquasivarieties:
\begin{prop}\label{P:5.4}
For any class $K$ of algebras, the following assertions hold:\\
(1) $H{\bf Q}(K) = {\bf S}{\bf P}{\bf P_{u}}{\bf D}(K)$;\\
(2) $H{\bf Q}(K) = {\bf S}{\bf P_{u}}{\bf P}{\bf D}(K)$;\\
(3) $H{\bf Q}(K) = {\bf S}{\bf P_{u}}{\bf P_{\omega}}{\bf D}(K)$;\\
(4) $H{\bf Q}(K) = {\bf S}{\bf L_{s}}{\bf P}{\bf D}(K) = {\bf L_{s}}{\bf S}{\bf P}{\bf D}(K) = {\bf L_{s}}{\bf P_{s}}{\bf S}{\bf D}(K)$.
\end{prop}

Via proposition 4.1 we obtain the following reformulation of Mal'cev
theorems of \cite[p. 214 - 215]{AIM}:
\begin{thm}\label{T:5.5}
A class $K$ of algebras of a given type is a hyperquasivariety if and only if
$K$\\
i) is ultraclosed;\\
ii) is heraditery;\\
iii) is multiplicatively closed; \\
iv) contains a trivial system; \\
v) is deriverably closed.
\end{thm}
\begin{thm}\label{T:5.6}
For every class $K$ of algebras of a given type we have
\begin{center}
$H{\bf Q}(K) = {\bf S}{\bf P_{r}}{\bf D}(K_{0})$
\end{center}
\end{thm}
\begin{rem}
If we accept, that the direct product of an empty family of algebras of a given
type  is a trivial algebra of a given type, we may remove condition iv) of
theorem  5.4 and substitute $K_{0}$ by $K$ in theorem 5.4.
(cf. Corollary 2.3. of \cite[p. 78]{VAG}).
\end{rem}

\end{document}